\documentclass[12pt,oneside,reqno,english]{amsart}
\usepackage{amsmath}
\usepackage{amssymb,amsbsy,amsthm}

\usepackage[foot]{amsaddr}

\usepackage{graphicx}
\usepackage{mathtools}
\usepackage{caption}
\usepackage{algorithm}
\usepackage[noend]{algpseudocode}
\usepackage{enumerate}
\usepackage{cases}
\usepackage[margin=1in]{geometry}
\usepackage{pdfsync}
\usepackage{hyperref}
\usepackage{wrapfig} 

\allowdisplaybreaks[1]
\numberwithin{equation}{section}

\DeclareMathOperator{\N}{\mathbb{N}} 
\DeclareMathOperator*{\argmax}{arg\,max} 

\newcommand{\p}{\mathbb{P}} 
\newcommand{\eps}{\varepsilon} 

\newtheorem{theorem}{Theorem}[section]
\newtheorem{lemma}[theorem]{Lemma}
\newtheorem{proposition}[theorem]{Proposition}

\newtheorem{corollary}[theorem]{Corollary}
\newtheorem{conjecture}[theorem]{Conjecture}
\newtheorem{question}[theorem]{Question}
\newtheorem{problem}[theorem]{Problem}
\newtheorem{definition}[theorem]{Definition}

\newtheorem{example}[theorem]{Example}

\title{Word length, bias and bijections in Penney's ante}

\author{Mathew Drexel$^{1}$} 

\author{Xuanshan Peng$^{1}$} 

\author{Jacob Richey$^{2} \dagger$}

\address{$^1$The University of British Columbia, Vancouver, Canada}
\address{$^2$The Alfr\'ed R\'enyi Institute of Mathematics, Budapest, Hungary}
\address{$\dagger$jfrichey001@gmail.com}


\begin{document}\maketitle\thispagestyle{empty}


\begin{abstract} 
Fix two words over the binary alphabet $\{0,1\}$, and generate iid Bernoulli$(p)$ bits until one of the words occurs in sequence. This setup, commonly known as Penney's ante, was popularized by Conway, who found (in unpublished work) a simple formula for the probability that a given word occurs first. We study win probabilities in Penney's ante from an analytic and combinatorial perspective, building on previous results for the case $p = \frac{1}{2}$ and words of the same length. For words of arbitrary lengths, our results bound how large the win probability can be for the longer word. When $p = \frac{1}{2}$ we characterize when a longer word can be statistically favorable, and for $p \neq \frac{1}{2}$ we present a conjecture describing the optimal pairs, which is supported by computer computations. Additionally, we find that Penney's ante often exhibits symmetry under the transformation $p \to 1-p$. We construct new explicit bijections that account for these symmetries, under conditions that can be easily verified by examining auto- and cross-correlations of the words. 
\end{abstract} 



\section{Introduction} \label{int}

Penney's ante refers to the following random `race': pick two finite binary words $v = v_1 \cdots v_m$ and $w = w_1 \cdots w_n$ and generate an iid Bernoulli$(p)$ sequence $(X_n)_{n \geq 1}$. The race ends when either $v$ or $w$ occurs in $X$: that is, at time $\min\{\tau_v,\tau_w\}$, where for any binary word $u = u_1 \cdots u_r, \tau_u$ is the stopping time 

\begin{equation} \tau_u = \min \{t \in \N: X_{t-r+1} = u_1, X_{t-r+2} = u_2,  \ldots X_t = u_r\}. \end{equation}

\medskip

If word $v$ occurred first, i.e. $\tau_v < \tau_w$, then $v$ `wins the race'. Assuming that neither $v$ nor $w$ is a subword of the other, the race is non-trivial in that both words win with positive probability. A basic study of the probability distributions $\tau_v$, which are approximately geometric, appears in Feller's classic probability text \cite{feller1991introduction}. Conway popularized the `race' version as an example of a non-transitive tournament (on the set of possible words), and gave an explicit formula for the odds that one word occurs before another. 
The formula, summarized in Theorem \ref{calc1}, arises from the so-called {\em{Abracadabra}} martingale construction, which is folklore among probabilists (see for example Section 17.3.2 of \cite{levin2017markov}). Conway's formula was generalized by Guibas/Odlyzsko and Li \cite{guibasodlyzko1981overlaps}\cite{li1980martingale}, who found linear systems describing the win probabilities and relevant generating functions for a race among an arbitrary (finite) family of words, and later to underlying sequences generated from an arbitrary finite-state Markov chain \cite{Pozdnyakov2009}\cite{gerberli1981hittingtimes}. 

Our aim in this work is to give a birds eye view of the combinatorial structure of Penney's ante through the lens of the win probability $\p_p(\tau_v < \tau_w)$, where $\p_p$ denotes iid Bernoulli$(p)$ measure. The overarching questions we have in mind are: for fixed $p \in (0,1)$ and pairs of words $(v,w)$ satisfying some natural conditions, can we describe the maximum or minimum win probability, and the optimizing pairs $(v,w)$? Under what conditions on the words $v, w$ does the win probability, viewed as a function of $p$, have some special property, like monotonicity or symmetry, and when can those symmetries be explained by a bijection? 

For the unbiased case, previous works have identified the optimal choice among words of the same length and described the optimizers \cite{phillips2021strategies}\cite{Csirik_1992}. Results are also available in the case where the underlying digit sequence is iid from a distribution with finite support \cite{Stark_1995} \cite{felix2006optimal}. We extend these results to words of different lengths, focusing on the question: when can a longer word possibly beat a shorter one? We prove that at $p = \frac{1}{2}$, a longer word can be a favorable choice only in a select few cases. Our findings are in agreement with a well-known rule of thumb for how to choose a word $v$ to beat a given word $w$ (see Theorem 1.6 of  \cite{guibasodlyzko1981overlaps}). In the case $p \neq \frac{1}{2}$, we give a conjecture along similar lines, which is supported by computer computations. Interestingly, in this case the optimum win probability appears to occur for words that are essentially periodic with period $1$ or $2$, depending on the value of $p$. We do not fully understand this phenomenon -- see Sections \ref{longlongbias} and \ref{open} for further discussion.

\begin{figure} \label{longshort_monoconj}
\includegraphics[width=.85\textwidth]{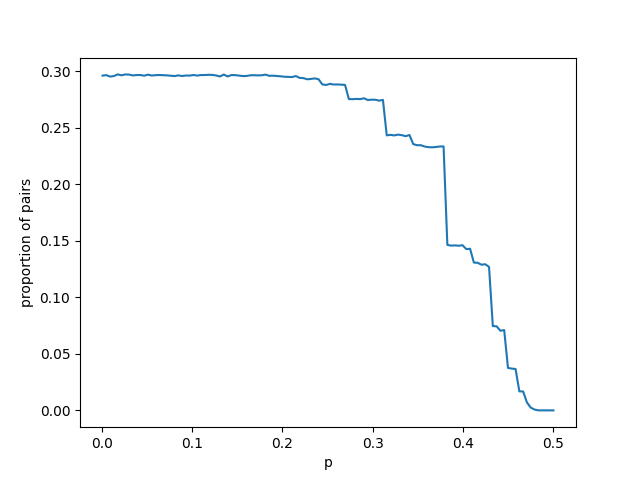} 
\centering
\captionsetup{margin=.8cm}
\caption{{\small{The approximate proportion of pairs of length at most $n = 17$ where the longer word is statistically favorable over the shorter word, plotted for $p \in (0, \frac{1}{2})$, as in Conjecture \ref{monoconj}.}}} \label{fig:monoconj}
\end{figure}

One naturally expects that if $p$ is perturbed away from $\frac{1}{2}$ -- decreased, say -- then it becomes a more common occurrence for a longer word to beat a shorter one. For $p$ close to 0, we show in Theorem \ref{bread} that the winner is essentially determined by which word has fewer 1s, disregarding the total length of the word (though among words with the same number of $1$s, there is a small subtlety). One way to formalize this heuristic is the following conjecture, which is supported by computer computations (see Figure \ref{fig:monoconj}):

\begin{conjecture} \label{monoconj} For sufficiently large integer $n$, the number of pairs of binary words $(v, w)$ of length at most $n$ such that $v$ is longer than $w$ and $\p_p(\tau_v < \tau_w) > \frac{1}{2}$ is non-increasing as a function of $p \in \left(0, \frac{1}{2}\right)$. \end{conjecture}

A natural guess for how to prove Conjecture \ref{monoconj} is to show that for each fixed pair $(v, w)$, the win probability $\p_p(\tau_v < \tau_w)$ is a monotone (or constant) function of $p \in \left(0, \frac{1}{2}\right)$. Surprisingly, this is false: a small example is the pair $(v,w) = (100010, 001100)$, with win probability $\p_p(\tau_v < \tau_w) = \frac{1}{2} \left(1-p(1-p)^4\right)$ taking its minimum value at $p = \frac{1}{5}$. See Section \ref{open} for further discussion. 

For certain pairs $(v, w)$, the win probability exhibits some symmetry under the transformation $p \to 1-p$, i.e. exchanging $0$s and $1$s. There are some obvious ways this can occur: as an example, if the word $w$ is obtained from $v$ by exchanging all $0$s to $1$s and vice versa, then exchanging all $0$s and $1$s in both the sequence $X$ and the words $v, w$ shows that $\p_p(\tau_v < \tau_w) = \p_{1-p}(\tau_w < \tau_v)$. In the language of probability, flipping all the bits gives a coupling between Penney's ante played between the pair $v, w$ and the pair $\overline{v}, \overline{w}$, the bit-flips of $v$ and $w$, such that $v$ wins the $v, w$ race if and only if $\overline{w}$ wins the $\overline{v}, \overline{w}$ race (where $\overline{u}$ is the bit-flip of $u$). 

%

Of course, using Conway's formula (Equation \ref{conway}), one can always express the win probability as a rational function of $p$, and explicitly check whether some symmetry holds. However, we find that many symmetries occur for non-obvious reasons; and crucially, since the proof of Conway's formula relies (at its root) on the optional stopping theorem for martingales, it does not naturally yield bijective proofs of these kinds of identities. This leads us to study the underlying combinatorial structure directly: namely, we consider the set $\Omega$ of all finite binary words that avoid the words $v$ and $w$, and seek explicit bijections between certain subsets of $\Omega$. Specifically, we will seek bijections $\varphi$ between sets of such words that either preserve the iid measure $\p_p$ or exchange $\p_p$ with $\p_{1-p}$, and are locally computable: that is, to determine a given digit of the an image word $\varphi(x)$, one should only have to check a bounded-sized window of $x$, or a window of finite expected size. This is similar to the approach taken by the authors of \cite{DEM21}, who show (by exact counting) a symmetry for the pair $(111, 101)$ and a specific value of $p$. In Proposition \ref{replace} and Theorem \ref{cake} we present large families of pairs for which we construct new, explicit, local bijections. We also find an explicit bijection for the pair $(v, w) = (1100, 1010)$ which arises from a natural graph construction, and demonstrates for any $p \in (0,1)$ the equality $\p_p(\tau_v < \tau_w) = \p_{1-p}(\tau_v < \tau_w)$ (both are equal to $\frac{1+p(1-p)}{2+p(1-p)}$). 

\medskip

{\bf{Remark}}: Our results are tangential to {\em{symbolic dynamics}}, a field at the intersection of combinatorics, dynamics, and ergodic theory -- see for example \cite{LM21}. The central objects in symbolic dynamics are {\em{shift spaces}}, akin to the sets $\Omega$, and the relevant isomorphisms are {\em{conjugacies}}, akin to the maps $\varphi$. We note that the condition appearing in Theorem \ref{replace} gives a class of simple conjugacies, while the more interesting bijections constructed in the proof of Lemma \ref{special} are not conjugacies (since they are not locally computable). Typically in symbolic dynamics one studies the {\em{measure(s) of maximal entropy}} for a shift space -- here we work with the measures $\p_p$, which are never maximal entropy measures for the corresponding shift spaces.

\subsection{Overview}

This work is structured as follows. Section \ref{bg} introduces the notation and formulas involved, most important among them being the Conway formula, Equation \ref{conway}, for the win probability. In Section \ref{res} we briefly advertise our main results and their implications. Section \ref{longlongman} covers the results and conjectures regarding words of different lengths, and Section \ref{lotto} presents bijections between different values of $p$ for fixed words of the same length. In Section \ref{open} we discuss some conjectures, questions, and possible future directions. 

\vspace{6pt}
\section{Preliminaries} \label{bg}


For a word $v \in \{0,1\}^n$, we often write $v$ in the shorthand $v = v_1 v_2 \cdots v_n$ for $v_i \in \{0,1\}$. For a letter $a \in \{0,1\}$, $\overline{a}$ denotes the bit-flip of $a$, $\overline{a} = 1-a$, and $||v||_a$ denotes the number of digits of $v$ that are equal to $a$. For a word $v$, $\overline{v}$ denotes the bit-flip of $v$, or the word $\overline{v} = \overline{v_1} \overline{v_2} \cdots \overline{v_n}$. For a pair of binary words $x, y$, we often write $xy$ for the concatenation of $x$ and $y$. For $a \in \{0,1\}$, we use the shorthand $a^n$ for the word consisting of $n$ copies of $a$. 

\begin{definition} \label{def:omegas} For binary words $v,w$, let $\Omega_{v,w}(v)$ be the set of all binary words $u = xv$ such that $v$ appears exactly once as a subword of $u$ (as the suffix). Define $\Omega_{w,v}(w)$ similarly, and set $\Omega_{v,w} = \Omega_{v,w}(v) \cup \Omega_{v,w}(w)$. Also, set $\Omega^*_{v,w}(v) = \{u: uv \in \Omega_{v,w}(v)\}$ (where $uv$ is the concatenation of words $u$ and $v$), and define $\Omega^*_{v,w}(w)$ similarly. \end{definition} 

$(\Omega_{v,w}, \p_p)$ is the natural probability space on which the random variables $\tau_v$ and $\tau_w$ are defined in the context of Penney's ante. 

\begin{definition} \label{def:overlappies} (Overlap sets and polynomials) Given words $v$ and $w$, the overlap set $\mathcal{O}(v,w)$ is the set of prefixes of word $w$ that are also suffixes of word $v$: 

\begin{equation} \mathcal{O}(v,w) = \{w_1 w_2 \cdots w_r: v_{m-r+1} v_{m-r+2} \cdots v_m = w_1 w_2 \cdots w_r\}. \end{equation}

For any $p \in (0,1)$, the correlation polynomial $vw_p$ is the real number 

\begin{equation} vw_p = \sum_{u \in \mathcal{O}(v,w)} \p_p(u)^{-1} = \sum_{u \in \mathcal{O}(v,w)} p^{-||u||_1} (1-p)^{-||u||_0}. \end{equation} \end{definition}

Simply put, $\mathcal{O}(v, w)$ is the set of words $u$ that are both a suffix of $v$ and a prefix of $w$. The overlap polynomial pops out of a clever martingale construction, sometimes called the {\em{Abracadabra}} martingale \cite{Sov66}\cite{li1980martingale}. Where it does not cause confusion we abuse notation and use $xy$ for concatenation of words, and $xy_p$ for the correlation polynomial value, since both are common in the literature.


Our main tool is the following formula, originally derived by Conway in unpublished work (see \cite{li1980martingale} for detailed proofs).  

\begin{theorem}[Win probability] \label{calc1} Let $v, w$ be any words over the binary alphabet $\{0,1\}$. Then 

\begin{equation} \mathbb{E}_p \tau_v = vv_p, \end{equation}

and if neither $v$ nor $w$ is a subword of the other,  

\begin{equation} \label{conway} {\normalfont{\text{Win}}}(v,w;p) \vcentcolon= \p_p(\tau_v < \tau_w) = \p_p(\Omega_{v,w}(v)) = \frac{ww_p - wv_p}{ww_p + vv_p - wv_p - vw_p}. \end{equation} \end{theorem}


It is not possible for arbitrary sets to occur as overlap sets. We will use the following special property of overlap sets, which is a simple exercise: 


\begin{theorem}[Forward propagation \cite{guibasodlyzko1981periods}] \label{fprop} Let $v$ be a binary word of length $n$. The set 

\begin{equation} \{n - i: v_1 v_2 \cdots v_i \in \mathcal{O}(v,v)\}\end{equation}

is closed under addition. \end{theorem}

Additionally, in Section \ref{lotto} we will use the following `bad prefix' sets: 

\begin{definition} \label{def:badprefixes} For words $v, w$, let $D(v,w)$ denote the set of non-empty prefixes of $v$ corresponding to non-trivial overlaps of $v,w$: that is, $D(v,w)$ consists of all words $v_1 \cdots v_{n-k}$ such that $w_1 \cdots w_k = v_{n-k+1} \cdots v_n$, i.e. $w_1 \cdots w_k \in \mathcal{O}(v, w)$. Also, let $F(v,w) \vcentcolon = D(v,v) \cup D(w,w) \cup D(w, v) \cup D(v,w)$.\end{definition} 

Note that $\mathcal{O}(v,w)$ and $D(v,w)$ are not symmetric in $(v,w)$, but $F(v,w)$ is. To clarify the definitions, we work out everything for an explicit example. 

\begin{example} Consider the pair $(v, w) = (1101110, 0110)$. We have
\begin{equation} \mathcal{O}(v,v) = \{1101110, 110\}, \hspace{5pt} \mathcal{O}(w,w) = \{0110, 0\}, \hspace{5pt} \mathcal{O}(v,w) = \{0\}, \hspace{5pt} \mathcal{O}(w, v) = \{110\}, \end{equation}
which gives 
\begin{equation} vv_p = p^{-5}(1-p)^{-2} + p^{-2}(1-p)^{-1}, \hspace{5pt} ww_p = p^{-2}(1-p)^{-2} + (1-p)^{-1}, \end{equation}
\begin{equation}\hspace{5pt} vw_p = (1-p)^{-1}, \hspace{5pt} wv_p = p^{-2}(1-p)^{-1}, \hspace{5pt}  {\normalfont{\text{Win}}}(v,w;p) = \frac{p^4(1+p-p^2)}{1+p^3}.\end{equation}

Similarly, 
\begin{equation} D(v,v) = \{1101\}, \hspace{5pt} D(w,w) = \{011\}, \hspace{5pt} D(v,w) = \{110111\}, \hspace{5pt} D(w, v) = \{0\}, \end{equation}

and $F(v,w) = \{0, 011, 1101, 110111\}.$ 
\end{example}


\vspace{6pt}

\section{Overview \& Results} \label{res}

Recall that for $p \in (0,1)$, $\p_p(\tau_v < \tau_w) = {\normalfont{\text{Win}}}\left(v,w;p\right)$ is the probability that $v$ beats $w$ in Penney's ante. Our first main result, proved in Section \ref{longlongman}, characterizes when a longer word can be statistically favorable over a shorter one for $p = \frac{1}{2}$. 

\begin{theorem} \label{tomatosoup} Let $v$ and $w$ be words over the alphabet $\{0, 1\}$ with len$(v) = $ len$(w) + 1$, and such that $v$ does not contain $w$ as a subword. If ${\normalfont{\text{Win}}}\left(v,w;\frac{1}{2}\right) \geq \frac{1}{2}$, then $w \in \{1^n, 0^n\}$. \end{theorem}

Since by Theorem \ref{calc1} the words $1^n$ and $0^n$ have the largest possible expected hitting time of all words of length $n$, it is intuitive that longer words can have the best chance against these two candidates. We also calculate the optimal win probability and the optimizing pair for a word that is $k$ digits longer than its opponent: 

\begin{theorem} \label{grilledcheese} Let $v$ and $w$ be such that len$(v) = $ len$(w) + k$ for non-negative integer $k$, and $v$ does not contain $w$ as a subword. Then 

\begin{equation} {\normalfont{\text{Win}}}\left(v,w;\frac{1}{2}\right) < \frac{2}{1+2^k} \end{equation}

and equality achieved asymptotically as $n \to \infty$ for the pair $(v, w) = (0^{k+1} 1^{n-1}, 1^n)$. \end{theorem} 

It follows immediately that for $p = \frac{1}{2}$ and $k \geq 2$, a word of length $n + k$ can never have favorable odds against a word of length $n$. A well-known rule appearing in previous work is that to optimize the win probability against a fixed word $w$, you should take a long prefix of $w$ and make it a suffix of $v$ \cite{guibasodlyzko1981overlaps}. Heuristically, this is a good choice because sometimes, when $w$ is about to occur, $v$ has just occurred. We note that the proof of Theorem \ref{tomatosoup} and the explicit pair appearing in Theorem \ref{grilledcheese} show that this rule of thumb holds for the case $k \geq 1$. 

The case $p \neq \frac{1}{2}$ is qualitatively different from $p = \frac{1}{2}$ in the following way. At $p = \frac{1}{2}$, the probability of a word depends only on its length, but this is false when $p \neq \frac{1}{2}$. Thus, words with higher probability under $\p_p$ will have an advantage in Penney's ante, regardless of their length. But if two words have comparable probabilities under $\p_p$, then the overlap patterns will be the key feature that determines the winner. This suggests restricting our attention to pairs $(v,w)$ whose probabilities are comparable under $\p_p$ for all $p$. One natural way to achieve this is to require $v$ and $w$ to have the same number of $1$s. Thus we seek to optimize over the following set: 

\begin{definition} For any integer $n \geq 2$ and $k \geq 1$, let $W(n,k)$ denote the set of pairs $(v, w)$ of words of length $n+k$ and $n$ respectively such that $||v||_1 = ||w||_1$ and $w$ is not a subword of $v$. \end{definition}

Computer calculations suggest that the optimal pairs in $W(n,k)$ can be identified as follows. We discuss this conjecture and its implications further in Section \ref{longlongbias}.  

\begin{conjecture} \label{baymax} For $k \in \{1, 2, \ldots\}$ and $n$ sufficiently large, there exists $r = r(n,k)$ such that: 

\begin{itemize} \item If $k = 1$, then 

\begin{equation} \argmax_{(v,w) \in W(n,k)} {\normalfont{\text{Win}}}(v,w;p) = \begin{cases} (001^{n-1}, 1^{n-1}0), & 0 < p \leq r(n,1) \\ (101^{n-1}, 1^n), & r(n,1) \leq p < 1 \end{cases} \end{equation}

\item If $k \geq 2$, then 

\begin{equation} \argmax_{(v,w) \in W(n,k)} {\normalfont{\text{Win}}}(v,w;p) = \begin{cases}  (0^{k+1} (10)^{m-1} 1, (10)^m), & 0 < p \leq r(n,k) \text{ and } n = 2m \\ (0^{k+1}(01)^{m}, (01)^m 0), & 0 < p \leq r(n,k) \text{ and } n = 2m+1 \\ (0^{k-1} 1 0 1^{n-1}, 1^n), & r(n,k) < p < \frac{1}{2} \\\end{cases} \end{equation}

where $m$ is an integer.

\end{itemize}

\end{conjecture}


We now turn to properties of the win probability, viewed as a function of $p \in (0,1)$. For a fixed pair $v, w$ of words with the same length, we show that the limit $p \to 0$ can only give three possible win probabilities: 

\begin{proposition} \label{bread} For any words $v, w$, $\lim_{p \to 0} {\normalfont{\text{Win}}}(v,w;p) \in \{0, \frac{1}{2}, 1\}.$ Moreover, the limit is $\frac{1}{2}$ if and only if $||v||_1 = ||w||_1 = r$ and neither $\mathcal{O}(v,w)$ nor $\mathcal{O}(w,v)$ contains an overlap word $u$ such that $||u||_1 = r$. \end{proposition}

The condition in the second part of Proposition \ref{bread} should be thought of as a weaker form of the condition that neither $v$ nor $w$ contains the other as a subword. In a nutshell, it says that for $p \approx 0$, $v$ and $w$ do not occur almost at the same time with high probability. 

Next, we focus on three possible symmetries of the function $\text{Win}(v,w;p)$. Say a pair $v,w$ has {\em{even}}, {\em{odd}}, or {\em{constant symmetry}} if $\text{Win}(v,w;\frac{1}{2} + x)$ is an even, odd, or constant function for $x \in (-\frac{1}{2}, \frac{1}{2})$, respectively. Note that by Proposition \ref{bread}, even (or constant) symmetry can only occur if Win$(v,w;p) = \frac{1}{2}$ for all $p \in (0,1)$ (if $v \neq w$). While these symmetries can occur for trivial reasons -- for example, $w = \overline{v}$ easily implies odd symmetry -- most of the time it is not obvious why we observe such a symmetry. Our main results, Theorems \ref{replace} and \ref{cake}, identify large classes of pairs where these symmetries occur, and give explicit bijections which explain the symmetry.

\begin{theorem} \label{replace} Let $\mathcal{R}$ denote the family of pairs of binary words $(v,w)$ of the same length such that $\Omega^*_{v,w}(v) = \Omega^*_{v,w}(w)$ (see definition \ref{def:omegas}).  

\begin{enumerate}[i.] \item For $(v,w) \in \mathcal{R}$,

\begin{equation} \label{eq:repl} {\normalfont{\text{Win}}}(v,w;p) = \frac{\p_p(v)}{\p_p(v) + \p_p(w)} = \frac{p^s}{p^s + (1-p)^s}, \end{equation} 

where $s = ||v||_1 - ||w||_1.$ 

\item If $(v,w) \in \mathcal{R}$, then $(v,w)$ has odd symmetry; and if in addition $||v||_1 = ||w||_1$, then $(v,w)$ has constant symmetry. 

%
%

\item The condition $(v,w) \in \mathcal{R}$ can be checked in polynomial time (with respect to the length of the words). 

\end{enumerate} \end{theorem} 

The condition $\Omega^*_{v,w}(v) = \Omega^*_{v,w}(w)$ says that the copy of $v$ (resp. $w$) that occurs right at the end of the string race can be replaced with $w$ (resp. $v$) without creating additional copies of $v$ or $w$. Thus in this case the bijection we seek is simply the identity map on $\Omega^*_{v,w}(v) = \Omega^*_{v,w}(w)$. A computer check suggests that property $\mathcal{R}$ is macroscopically common among all pairs of words: 

\begin{conjecture} \label{ratio1} Let $\mathcal{R}_n$ denote the sub-family of pairs in $\mathcal{R}$ of length $n$. The limit

\begin{equation} \lim_{n \to \infty} 4^{-n} |\mathcal{R}_n| \end{equation} 

exists and is strictly positive. 

\end{conjecture} 

It is easy to show that if $v$ and $w$ have trivial self and cross overlap sets, i.e. if $\mathcal{O}(v,w) = \mathcal{O}(w,v) = \emptyset$ and $\mathcal{O}(w,w) = \mathcal{O}(v,v) = \{n\}$, then $(v, w) \in \mathcal{R}$, and Equation \ref{conway} gives the conclusions $i$ and $ii$. However, an exact computer search with $n = 15$ returned the approximate value $.15753$ for the limit in Conjecture \ref{ratio1}, while the proportion of pairs that have trivial auto- and cross- overlap sets is approximately $.02309$. Thus, Theorem \ref{replace} gives a much stronger result than than one obtains by a naive application of Equation \ref{conway}. The pair $v =000100, w=001110$ is an example where property $\mathcal{R}$ holds -- see Lemma \ref{repcheck} -- but the auto- and cross-correlations are non-trivial.

Still, property $\mathcal{R}$ is not necessary to have odd symmetry. The minimal example is $1000$ and $0110$: this pair has odd symmetry, with Win$(1000, 0110;p) = 1-p$ -- which is what one would obtain by applying Equation \ref{eq:repl} to $(v,w)$ -- but  $(v,w) \notin \mathcal{R}$. Oddly, $(v, \overline{w}) \in \mathcal{R}$, but this seems to be a coincidence. We do not know of a simple bijective argument that explains the odd symmetry in this case. See Section \ref{open} for further discussion.

Similarly, property $\mathcal{R}$ is not necessary for constant symmetry. The minimal example $(v,w) \notin \mathcal{R}$ with constant symmetry is $(01100101, 01010110)$, which we discovered in a computer search. We explain the symmetry for this example by an explicit, novel bijection in Section \ref{lotto}, and generalize that construction to obtain a large family of such pairs in the following Theorem. 

\begin{theorem} \label{cake} There exists a family $\mathcal{E}$ of pairs of words $(v,w)$ of the same length such that the following hold:

\begin{enumerate}[i.] \item For $(v,w) \in \mathcal{E}$, there exists a bijection $\varphi: \Omega_{v,w}(v) \to \Omega_{v,w}(w)$ such that $\p_p(x) = \p_p(\varphi(x))$. In particular, $(v,w)$ has constant symmetry. 

\item The condition $(v,w) \in \mathcal{E}$ can be checked in polynomial time (with respect to the length of the words). 

\item $\mathcal{E} \setminus \mathcal{R}$ is infinite. \end{enumerate} \end{theorem}

{\bf{Remark:}} We note (but do not prove) that the map $\varphi$ has a {\em{finite expected coding radius,}} meaning: given a random $X \in \Omega_{v,w}(v)$ distributed according to the measure $\mathbb{P}_p$ conditioned on the event $\tau_v < \tau_w$, any digit of the image word $\varphi(x)_i$ can be determined by examining a block $x_j x_{j+1} \cdots X_{j+R-1}$ of random size $R$ with $\mathbb{E} R < \infty$ (where the expectation is with respect to the conditioned measure). In the language of symbolic dynamics, this is the next best thing after being a conjugacy-like map, i.e. if $R$ were finite almost surely.

\medskip

Finally, we turn to even symmetry, which we understand the least. By using a graphical representation of Penney's ante, we construct an explicit bijection for the pair $(1010, 1100)$, which has even symmetry. Although similar bijections likely work to explain even symmetry for other pairs, we do not know how to generalize this bijection to a larger family of pairs -- see Section \ref{open} for further discussion. 

%
%
%
%
%

\vspace{6pt}
\section{Long vs short} \label{longlongman}

Our starting point is the proof of Theorem \ref{tomatosoup}, which characterizes when longer words can be favorable over shorter ones in the unbiased case $p = \frac{1}{2}$. 


\begin{proof}[Proof of Theorem \ref{tomatosoup}] We prove the contrapositive statement, so assume $w \neq 1^n$ or $0^n$. Recall the correlation polynomials (Definition \ref{def:overlappies}), and for simplicity, write $xy = xy_{\frac{1}{2}}$ for the auto- or cross correlation polynomial, with $x,y \in \{v,w\}$. By Theorem \ref{calc1} and a bit of simple algebra, 

\begin{equation} \label{tomatoeq1} \text{Win}(v,w;p) < \frac{1}{2} \iff ww - wv < vv-vw. \end{equation}

Note that $n-1 \notin \mathcal{O}(w,w)$, since otherwise by Theorem \ref{fprop}, $\mathcal{O}(w,w) = \{1, 2, \ldots, n-1, n\}$, which is easily seen to be equivalent to $w = 1^n$ or $0^n$.  

We now divide into a few simple cases, and establish equation \ref{tomatoeq1} in each one.

\begin{itemize} \item[\textbf{Case 1}] Suppose $n-1 \notin \mathcal{O}(v,w)$. Then  

\begin{equation} vw \leq \sum_{1 \leq i \leq n-2} 2^i = \frac{1}{2} 2^{n}. \end{equation}

If $n-2 \in \mathcal{O}(w,w)$, then $n-x \notin \mathcal{O}(w,w)$ for all odd integers $x$, since otherwise by Theorem \ref{fprop}, we would have $w_1 = w_3 = \cdots$, $w_2 = w_4 = \cdots$, and $w_1 = w_{x+1}$, which altogether imply $w = 1^n$ or $0^n$. Thus we have 

\begin{equation} ww \leq \left(\sum_{\substack{1\leq i \leq n \\ i \text{ even} }} 2^i \right) \wedge \left(\sum_{\substack{1 \leq i \leq n \\ i \neq n-2, n-1}} 2^i \right) \leq \frac{4}{3} 2^n. \end{equation}

We always have the bounds $vv \geq 2^{n+1}$ and $wv \geq 0$, so in this case we obtain 

\begin{equation} ww - wv \leq \frac{4}{3} 2^n \text{ and } vv - vw \geq 2^{n+1} - \frac{1}{2} 2^n = \frac{3}{2} 2^n.\end{equation} 

Since $\frac{3}{2} > \frac{4}{3}$, this gives the desired bound. 

\item[\textbf{Case 2}] Suppose $n-1, n-2 \in \mathcal{O}(v,w)$. Unraveling the definition of $\mathcal{O}(v,w)$ gives

\begin{equation} v = v_1 v_2 w_1 w_2 \cdots w_{n-1} = v_1 v_2 v_3 w_1 w_2 \cdots w_{n-2} \text{ and thus } w = w_1 w_1 \cdots w_1 w_n. \end{equation}

Since $w \neq 1^n$ or $0^n$, $w_n \neq w_1$, so $ww = 2^n$. Since $v$ does not contain $w$ as a subword, $n \notin \mathcal{O}(v,w)$, so $vw \leq \sum_{i=1}^{n-1} 2^i \leq 2^n$. Therefore 

\begin{equation} ww - wv \leq 2^n - 0 = 2^n, \text{ and } vv - vw \geq 2^{n+1} - 2^n = 2^n,\end{equation}

as desired.

\item[\textbf{Case 3}] Suppose $n-1 \in \mathcal{O}(v,w)$, $n-2 \notin \mathcal{O}(v,w)$, and $n-2 \notin \mathcal{O}(w,w)$. Then we bound $ww$ and $vw$ directly: 

\begin{equation} vw \leq \sum_{\substack{1 \leq i \leq n-1 \\ i \neq n-2}} 2^i \leq \frac{3}{4} 2^n, \end{equation} 

and 
\begin{equation} ww \leq \sum_{\substack{1 \leq i \leq n \\ i \neq n-2, n-1}} 2^i \leq 2^n + 2^{n-2} = \frac{5}{4} 2^{n}. \end{equation}

Combining these as before, with $vv \geq 2^{n+1}$ and $wv \geq 0$ gives the desired bound. 

\item[\textbf{Case 4}] Suppose $n-1 \in \mathcal{O}(v,w)$, $n-2 \notin \mathcal{O}(v,w)$, and $n-2 \in \mathcal{O}(w,w).$ Then $w_1 = w_3 = \cdots$ and $w_2 = w_4 = \cdots$, so without loss of generality,

\begin{equation} w = \begin{cases} (10)^{\frac{n}{2}}, & n \text{ even } \\ (10)^{\frac{n-1}{2}}1, & n \text{ odd} \end{cases} \end{equation}

and 

\begin{equation} v = \begin{cases} v_1 v_2 (10)^{\frac{n}{2} - 1}1, & n \text{ even } \\ v_1 v_2 (10)^{\frac{n-1}{2}}, & n \text { odd} \end{cases} \end{equation} 

In both cases, by a direct computation,  

\begin{equation} vw = \sum_{\substack{ 1 \leq i \leq n-2 \\ i \text{ odd } }} 2^{n-i} \leq \frac{2}{3} 2^n \end{equation} 

and 

\begin{equation} ww = \sum_{\substack{ 0 \leq i \leq n-1 \\ i \text{ even } }} 2^{n-i} \leq \frac{4}{3} 2^n. \end{equation}

Combining these as before, with $vv \geq 2^{n+1}$ and $wv \geq 0$ gives the desired bound.

\end{itemize} \end{proof}

Theorem \ref{tomatosoup} is tight, in the sense that the inequality Win $ \geq \frac{1}{2}$ cannot be replaced by \linebreak Win $ \geq \frac{1}{2} - \epsilon$ for any $\epsilon > 0$, since (for example) the pair $v = 11(10)^m$, $w = (10)^m 1$ achieves $q_{w,v}\left(\frac{1}{2}\right) = \frac{1}{2} - O(2^{-2m}).$ (This is the same pair that appears in \textbf{Case 4} of the proof of Theorem \ref{tomatosoup}.)




We now prove Theorem \ref{grilledcheese}, which gives a tight upper bound for the win probability in terms of the length difference $k$ between $v$ and $w$.  

\begin{proof}[Proof of Theorem \ref{grilledcheese}] We only need crude bounds on each of the overlap polynomials $vv, ww,$ $wv, vw$: namely, $vv \geq 2^{n+k}$, $ww \leq \sum_{i=1}^n 2^i \leq 2^{n+1}$, $wv \geq 0$, and since $v$ is not a subword of $w$, $n \notin \mathcal{O}(v,w)$, so $vw \leq \sum_{i=1}^{n-1} 2^i \leq 2^n$. Thus, by Theorem \ref{calc1},  

\begin{equation} \text{Win}\left(v,w;\frac{1}{2}\right) \leq \frac{2^{n+1}}{2^{n+k} + 2^{n+1} - 2^{n}} = \frac{2}{1 + 2^k}. \end{equation}

The pair $v = 0^{k+1} 1^{n-1}$, $w = 1^n$ asymptotically achieves each of the bounds on $vv, ww, wv$ and $vw$ up to negligible $o(1)$ terms.

\end{proof}


\subsection{Long vs short with bias} \label{longlongbias} 

We verified conjecture \ref{baymax} by an exact computer search for all pairs of words with $n \leq 11$ and $k \leq 4$. Note that all the optimal pairs appearing there have the property that $v$ has $w_1 w_2 \cdots w_{n-1}$ as a suffix, following the rule of thumb for the optimal choices in the case $p = \frac{1}{2}$ and $k = 0$. It is surprising that only words that are periodic (or contain a large, periodic subword) with period $1$ or $2$ appear as the optimizers. A back of the envelope calculation suggests that pairs where $w$ has a longer period do not give higher win probabilities, but we could not find a proof along the same lines as the proof of Theorem \ref{grilledcheese}. Assuming the conjecture, we can precisely compute the maximum win probability for longer words, and the transition point between (pseudo-)period 1 and 2 optimizers: 

\begin{corollary} \label{nanobots} Assume Conjecture \ref{baymax} holds, and recall the notation there. Fix a pair of binary words $(v,w) \in W(n,k)$. 

\begin{itemize} \item For $k = 1$, 

\begin{equation} {\normalfont{\text{Win}}}(v,w;p) \leq \begin{cases} 1-p, & 0<p\leq r(n,1) \\ \frac{1}{2-p}, & r(n,1) \leq p < \frac{1}{2} \end{cases} \end{equation}

The threshold satisfies

\begin{equation} r(n,1) = \kappa + O(\kappa^{n}) \text{ as } n \to \infty, \text{ where } \kappa = \frac{3-\sqrt{5}}{2}.\end{equation}

\item For $k \geq 2$, 

\begin{equation} {\normalfont{\text{Win}}}(v,w;p) \leq \begin{cases} \frac{(1-p)^k}{1-p(1-(1-p)^k) + p^2}, & 0 < p \leq r(n,k) \\ \frac{(1-p)^{k-1}}{1 + (1-p)^k}, & r(n,k) \leq p < \frac{1}{2} \end{cases} \end{equation}
\medskip
The threshold $r(n,k)$ is, up to $o(1)$ terms as $n \to \infty$, the unique root in $(0, \frac{1}{2})$ of the polynomial $(1-2z)(1-z)^k - z^2$. In particular, 

\begin{equation} r(n,k) \approx \frac{k+2-\sqrt{k^2-4k+8}}{4k-2} \text{ as } n, k \to \infty. \end{equation} 

\end{itemize}

\end{corollary}

We omit the proofs, which are straightforward calculations applying Theorem \ref{calc1} to the explicit pairs appearing in Conjecture \ref{baymax}. The upper bounds in Corollary \ref{nanobots} are continuous at the transition point $r(n,k)$. For the case where $k \geq 2$ and $n$ is odd, the optimal pair $(0^{k+1}(01)^{m}, (01)^m 0)$ has win probability which strictly increases as $n \to \infty$ to the value for the $n$ even case. 

\vspace{6pt}
\section{Symmetry and bijections} \label{lotto}

In this section we investigate symmetries that underlie the combinatorics of Penney's ante. We restrict our attention to words $v$ and $w$ with the same length. The most general question we could seek to answer is: for two values $p, p' \in (0,1)$ and for two pairs $(v,w), (v', w')$, when is there a coupling between the $v$ vs $w$ race under $\p_p$ and the $(v', w')$ race under $\p_{p'}$? Here we restrict our attention to comparing a fixed pair $(v,w)$ to itself, but under different values of $p$. Our primary aim is to characterize when symmetries of the function $\text{Win}(v,w;p)$ occur, and to construct explicit bijections that explain the symmetry. Recall the terminology from Section \ref{res}: the pair $(v,w)$ has {\em{odd symmetry}} if $\text{Win}(v,w;1-p) = 1-\text{Win}(v,w;p)$ for $p \in (0,1)$; {\em{even symmetry}}, if $\text{Win}(v,w;1-p) = \text{Win}(v,w;p)$ for $p \in (0,1)$; and {\em{constant symmetry}} if $\text{Win}(v,w;p) = \frac{1}{2}$ for all $p \in (0,1)$. 


As a warm-up, we prove Proposition \ref{bread}. 




\begin{proof}[Proof of Proposition \ref{bread}] We start with the first conclusion. The result is clear if either $v$ or $w$ contains the other as a subword, so assume this is not the case. Let $r =\max\{||v||_1, ||w||_1\}$. Observe that a term of the form $p^{-||w||_1} (1-p)^{-j}$ for non-negative integer $j$ occurs exactly once in the overlap polynomial $ww_p$, when $j = ||w||_0$. Similarly, there is at most one term of the form $p^{-||v||_1} (1-p)^{-k}$ in $vv_p$ for non-negative integer $k$. Also, at most one of $vw_p$ or $wv_p$ may contain a term of the form $p^{-r}(1-p)^{-j}$ for non-negative integer $j$: indeed, if (without loss of generality) $vw_p$ contains a term of the form $p^{-r} (1-p)^{-j}$, say $u \in \mathcal{O}(v, w)$ with $||u||_1 = r$, then we must have $v = 0^a 1 u' 0^b$ and $w = 0^c 1 u' 0^d$ for some non-negative integers $a, b, c, d$ with $c \leq a$, and where $u'$ is a subword of $u$. Since neither $v$ nor $w$ contains the other as a subword, $d > b$, which implies that $\mathcal{O}(w, v)$ does not have an overlap with $r$ 1s. 

Using formula \ref{calc1}, the win probability is a rational function of the form 

\begin{equation} \text{Win}(v,w;p) = \frac{\sum_{0 \leq i, j \leq m} (\alpha_{ij} - \beta_{ij}) p^{-i} (1-p)^{-j}}{\sum_{0 \leq i, j \leq m} (\alpha_{ij} -\beta_{ij} + \gamma_{ij} -\delta_{ij}) p^{-i} (1-p)^{-j}}, \end{equation}

where the $\alpha, \beta, \gamma, \delta$ coefficients are $\{0,1\}$ valued. Let $r'$ denote the maximum power of $p^{-1}$ that appears in either sum (after any cancellations). By the observations in the previous paragraph, $r' = r$, the coefficient of $p^{-r}$ in the numerator is $0$ or $1$, and the coefficient of $p^{-r}$ in the denominator is either $1$ or $2$ (up to some bounded number of factors of $1-p$, which are $1+O(p)$ as $p \to 0$). Thus the only possible $p \to 0$ limits are the quotients $\{\frac{0}{1}, \frac{0}{2}, \frac{1}{1}, \frac{1}{2}\} = \{0, \frac{1}{2}, 1\}$. 

For the second part, if $||v||_1 \neq ||w||_1$, then the denominator has leading coefficient $1$ and the numerator has leading coefficient either $0$ or $1$, so the limit must be $0$ or $1$. Similarly, if $||v||_1 = ||w||_1 = r$ and one of $\mathcal{O}(v, w)$ or $\mathcal{O}(w, v)$ has an overlap that contains $r$ 1s, then the leading coefficient in the numerator is $0$ or $1$, and the leading coefficient in the denominator is $1$, so the limit cannot be $\frac{1}{2}$. 

\end{proof}

\subsection{Odd Symmetry} \label{sec:odd}

Recall that we say {\em{odd symmetry}} holds if 
\begin{equation} \text{Win}(v,w;p) = 1-\text{Win}(v,w;1-p) \text{ for all } p \in (0,1),\end{equation}

or equivalently, if 
\begin{equation} \p_p(\tau_v < \tau_w) = 1-\p_{1-p}(\tau_v < \tau_w) = \p_{1-p}(\tau_w < \tau_v) = \p_{p}(\tau_{\overline{w}} < \tau_{\overline{v}}). \end{equation}

Our aim is to construct canonical bijections that give rise to this equality. A class of pairs that obviously has odd symmetry is when $v = \overline{w}$. This odd symmetry is naturally realized via the bijection $\varphi: \Omega_{v, \overline{v}}(v) \to \Omega_{v, \overline{v}}(\overline{v})$ given by $\varphi(u) = \overline{u}$. Since $\p_p(\varphi(u)) = \p_{1-p}(u)$, this bijection induces the equality $\p_p(\tau_v < \tau_{\overline{v}}) = \p_{1-p}(\tau_{\overline{v}} = \tau_{v})$. Next we turn to a much larger class of pairs with odd symmetry.

\begin{definition}[Property $\mathcal{R}$] We say a pair of binary words $(v,w)$ has property $\mathcal{R}$, and write $(v,w) \in \mathcal{R}$, if $\Omega^*_{v,w}(v) = \Omega^*_{v,w}(w)$. \end{definition} 

Simply put, property $\mathcal{R}$ says that at the end of the Penney's ante race, swapping the final $v$ to a $w$ or vice versa does not create any additional copies of $v$ or $w$. To check whether property $\mathcal{R}$ holds, it is not necessary to directly check if the infinite sets $\Omega^*_{v,w}(v)$ and $\Omega^*_{v,w}(w)$ are equal. Instead, one only needs to examine the `bad prefix' set $F$: 

\begin{lemma} \label{repcheck} A pair of words $(v,w)$ has property $\mathcal{R}$ if and only if $fv, fw \notin \Omega_{v,w}$ for all $f \in F(v,w)$ (see Definition \ref{def:badprefixes}). \end{lemma}


\begin{proof} Let $n$ be the length of $v, w$. We will justify the following chain of equivalences:

\begin{align} \label{equiv1} \Omega^*_{v,w}(v) = \Omega^*_{v,w}(w) & \iff \text{ for all } x \in \Omega^*_{v,w}, (xv \notin \Omega_{v,w} \iff xw \notin \Omega_{v,w}) \\
& \label{equiv2} \iff \text{ for all } x \in \{0,1\}^{n-1}, (xv \notin \Omega_{v,w} \iff xw \notin \Omega_{v,w})\\
& \label{equiv3} \iff \text{ for all } f \in F(v,w), (fv \notin \Omega_{v,w} \iff fw \notin \Omega_{v,w}) \\
& \label{equiv4} \iff \text{ for all } f \in F(v,w), (fv, fw \notin \Omega_{v,w}) \end{align}

The equivalence \ref{equiv1} is clear. For \ref{equiv2}, the $\implies$ implication is clear since the latter condition is weaker than the former. For the reverse implication, we used the fact that since $x \in \Omega^*_{v,w}$ does not have $v$ or $w$ as a subword, no $v$ or $w$ can occur in $xv$ or $xw$ without intersecting the final $n-1$ digits of $x$, along with the fact that $\Omega^*_{v,w}$ contains all words of length $n-1$. For \ref{equiv3}, the $\implies$ implication is again clear since the latter condition is weaker than the former. For the reverse implication, note that if for example $xv \notin \Omega_{v,w}$, then $xv$ contains either a copy of $w$ or two copies of $v$, which implies that some suffix $f$ of $x$ is an element of $D(v,v) \cup D(w, v)$, and similarly for $xw$. The final equivalence follows from the fact that $fv, fw \in \Omega_{v,w}$ for any $f \in \{0,1\}^{n-1} \setminus F(v,w)$. 
\end{proof}

\begin{proof}[(Proof of Theorem \ref{replace})]

By property $\mathcal{R}$, let $\Omega^* = \Omega^*_{v,w}(v) = \Omega^*_{v,w}(w)$. By the definitions of Win$(v,w;p)$ and $\Omega^*_{v,w}(v)$ and using that $\p_p$ is iid,
\begin{equation} \text{Win}(v,w;p) = \p_p(\Omega_{v,w}(v)) = \p_p(\Omega^*) \p_p(v), \end{equation} 

and similarly, 
\begin{equation} \p_p(\Omega_{v,w}(w)) = \p_p(\Omega^*) \p_p(w). \end{equation}

Since $v \neq w$, $\p_p(\Omega_{v,w}(v)) + \p_p(\Omega_{v,w}(w)) = 1$, and thus 

\begin{equation} \text{Win}(v,w;p) = \frac{\p_p(\Omega^*)\p_p(v)}{\p_p(\Omega^*)\p_p(v) + \p_p(\Omega^*)\p_p(w)} = \frac{\p_p(v)}{\p_p(v) + \p_p(w)}. \end{equation}

\medskip

Straightforward algebra gives the last equality of $i$, which directly implies odd symmetry, and constant symmetry if $\p_p(v) = \p_p(w)$, and thus $ii$ holds. Since the set $F(v,w)$ can be determined by checking equality of linearly many overlaps with respect to the length of $v$ ($=$ length of $w$), and similarly for a condition of the form $fv \in \Omega_{v,w}$ or $fw \in \Omega_{v,w}$ where $f$ has length at most that of $v$, part $iii$ follows from Lemma \ref{repcheck}. 

\end{proof}

\subsection{Constant symmetry} \label{sec:const} Recall that we say a pair $v,w$ has {\em{constant symmetry}} if \newline ${\normalfont{\text{Win}}}(v,w;p) = \frac{1}{2}$ for all $p \in (0,1)$. Note that Theorem \ref{replace} guarantees a large class of pairs with constant symmetry: $(v,w)$ has constant symmetry whenever $||v||_1 = ||w||_1$ and property $\mathcal{R}$ holds for the pair $(v,w)$. (An example is $0101, 0110$.) However, property $\mathcal{R}$ is not necessary for constant symmetry. We start by detailing the minimal such example, which we will then generalize to a large family of such examples in Theorem \ref{cake}. 

\begin{lemma} \label{special} Let $a = 0110, b = 0101$, and $v = 01100101 = ab, w = 01010110 = ba$. There exists a bijection $\varphi: \Omega_{v,w}(v) \to \Omega_{v,w}(w)$ such that for any $p \in (0,1)$, $\p_p(x) = \p_p(\varphi(x))$. The pair $(v,w)$ has constant symmetry.
\end{lemma}

Note that property $\mathcal{R}$ does not hold for $v,w$, since $a \in \Omega^*_{v,w}(v), \text{ but } a \notin \Omega^*_{v,w}(w)$. The rough idea is: we want to construct a bijection that swaps occurrences of $a$ and $b$. Although there are many candidate maps that seem to achieve this, we only found one bijection that works. 

\begin{proof} Note that each word $z \in \Omega_{v,w}(v)$ can be uniquely expressed in the form $z = x a^j v \in \Omega_{v,w}(v)$ where $j \geq 0$ is maximal. Write $j = j(z)$ for the maximal such $j$. Then define the map $\varphi_{a,b}$ by 
\begin{equation} \varphi_{a,b}(x a^j v) = x b^j w,  \end{equation}

and define $\varphi_{b,a}$ similarly, with the roles of $a$ and $b$ reversed. Note that $\p_p(a) = \p_p(b)$, and thus $\p_p(v) = \p_p(w)$, so both $\varphi_{a,b}$ and $\varphi_{b,a}$ are measure-preserving. Also, for $z = x a^{j(z)} v \in \Omega_{v,w}(v)$, $x$ cannot not have $b$ as a suffix (else $z$ would contain a copy of $w$), and thus $j(\varphi_{a,b}(z)) = j(x b^{j(z)} w) = j(z)$. (And the corresponding fact holds for the reverse map.) Then, assuming that $\varphi_{a,b}$ maps $\Omega_{v,w}(v)$ to $\Omega_{v,w}(w)$ (and vice versa for $\varphi_{b,a}$) -- this is the crux of the argument and will be proved below -- we would have that $\varphi_{a,b} = \varphi_{b,a}^{-1}$, so $\varphi_{a,b}$ is a measure-preserving bijection, completing the proof of the first part of the Lemma. We will show that $\varphi_{a,b}$ has image in $\Omega_{v,w}(w)$; the details for $\varphi_{b,a}$ are essentially the same. 

Write $\Omega_{v,w} = \Omega$ for short, and for a word $u$, let $u_{[i]} = u_1 u_2 \cdots u_i$ denote the $i$ prefix of $u$. We must show that if $yv \in \Omega(v)$, then $\varphi(yv) \in \Omega(w)$. The key observation is that for any word $y$ that does not contain $v$ or $w$ as a subword, 
\begin{equation}\label{sufeq1} yv \in \Omega(v) \iff y \neq y' u \text{ for } u \in \mathcal{D}(v, v) \cup \mathcal{D}(w,v) = \{b, v_{[6]}, w_{[7]}\} \end{equation}

and 
\begin{equation} \label{sufeq2} yw \in \Omega(w) \iff y \neq y' u \text{ for } u \in \mathcal{D}(w, w) \cup \mathcal{D}(v, w) = \{a, v_{[6]}, w_{[7]} \}. \end{equation}

Consider a word $z = xa^j v \in \Omega(v)$ with $j(z) = 0$, i.e. a word $x v \in \Omega(v)$ where $x$ does not end in $a$. By statement \ref{sufeq1}, $x$ does not end in $a, b, v_{[6]}$ or $w_{[7]}$; and then by statement \ref{sufeq2}, $xw \in \Omega(w)$ and $x$ does not end in $b$. So $\varphi_{a,b}$ is well-defined on the set of $x v \in \Omega(v)$ with $j(xv) = 0$.

Now consider a word $z = x a^j v \in \Omega(v)$ with $j(z) \geq 1$. Since $j$ is maximal, by statement \ref{sufeq1}, $x$ does not end in $a, b, v_{[6]}$ or $w_{[7]}$. Observe that 
\begin{equation} \mathcal{D}(v, bb) = \{a, v_{[6]}\} \text{ and } \mathcal{D}(w, bb) = w_{[7]}. \end{equation}

Since we assumed $x$ doesn't end in any of these words, we obtain that $xb^jw = xb^{j+1} a$ doesn't contain any copy of $v$ or $w$ that includes digits of $x$. Finally, a direct check of the words $b^3$ and $bba = bw$ shows that neither $v$ nor $w$ can occur as a subword of $b^j$ or $bw$. Combining all this shows $x b^j w \in \Omega(w)$, as desired. 

Since $\p_p(v) = \p_p(w)$, the existence of $\varphi_{a,b}$ implies (by the same calculation as in the proof of Proposition \ref{replace}) that $\text{Win}(v,w;p) = \frac{1}{2}$. \end{proof}

%

%

The construction in Lemma \ref{special} generalizes to the following class of words. 

\begin{definition}[Property $\mathcal{E}$] We say a pair of binary words $(v,w)$ has property $\mathcal{E}$, and write $(v,w) \in \mathcal{E}$, if there exist binary words $a, b$ and a finite set $S$ of binary words such that the following hold:

\begin{enumerate}[\normalfont I.] \item $||a||_1 = ||b||_1$
\item $v = ab, w = ba$
\item $D(v,v) \cup D(w,v) = S \cup \{a\} $ and $D(w,w) \cup D(v,w) = S \cup \{b\}$
\item For all $m \geq 2$, $D(v,a^m) \cup D(w,a^m) \subset S \cup \{a, b\}$ and $D(v, a^m) \cup D(w, b^m) \subset S \cup \{a, b\}$
\item For all $m \geq 2$, $a^m$ and $b^m$ do not contain $v$ or $w$ as a subword
\item For all $m \geq 1$, $a^m b \in \Omega_{v,w}(v)$, $b^m a \in \Omega_{v,w}(w)$. 
\end{enumerate} \end{definition}

\begin{proof}[(Proof of Theorem \ref{cake})] The argument in Lemma \ref{special} adapts immediately to this setting: properties II-VI imply, by the same argument, that the map $\varphi_{a,b}: \Omega_{v,w}(v) \to \Omega_{w,v}(w)$ given by
\begin{equation} 
\varphi_{a,b}(xv) = x' b^j w \text{ where } j = \max \{j': a^{j'} \text{ is a suffix of } x\}
\end{equation}

\medskip
is a well-defined bijection, and property I guarantees it is measure preserving. To see that the condition $(v,w) \in \mathcal{E}$ can be checked in polynomial time in the lengths of $v$ and $w$, first note that the existence of $a, b$ satisfying II can be determined by checking all prefixes of $v$, and if some such pair exists, it will be determined by the same algorithm. Given $a$ and $b$ satisfying condition II, conditions I and III require checking linearly many overlap comparisons, and if III holds, the same algorithm will determine the set $S$. Given $a, b$ and the set $S$, note that it suffices to check conditions IV, V and VI for $m < M$, where $M$ is at most a linear function of the lengths of $a$ and $b$. For example, the condition $D(v, a^m) \subset S \cup \{a, b\}$ does not depend on $m$ if $a^m$ has length greater than that of $v$, and similar facts hold for all the conditions in IV, V and VI. Thus at most polynomially many comparisons are required. Finally, an easy computation shows that for all $k \geq 2$, words $a = 01^k0, b = 01^{k-1}01$ give rise to a pair $(a_k b_k, b_k a_k) \in \mathcal{E} \setminus \mathcal{R}$. \end{proof}


\medskip

\subsection{Even symmetry} \label{sec:even}

We do not know a general construction akin to the families $\mathcal{R}$ or $\mathcal{E}$ which exhibit even (but not constant) symmetry. We study one example in depth, and give an explicit bijection in terms of a graphical representation of the space $\Omega_{v,w}$.

\begin{proposition} \label{prop:quad} Let $v = 1100, w = 1010$. There exists a bijection $\varphi: \Omega_{v,w}(v) \to \Omega_{v,w}(v)$ such that $\p_p(\varphi(x)) = \p_{1-p}(x)$ for all $p \in (0,1)$. In particular, $(v,w)$ has even symmetry.


\end{proposition}

We note that the existence of such a bijection does not come for free from the even property, or from the Conway formula. To prove Proposition \ref{prop:quad}, we rely on a certain graph representation of the string race, which appears in Figure \ref{1100_1010_graph}. This construction is a common way to compute recursions, probabilities and expectations explicitly in Penney's ante \cite{nickerson2007penney}. (In the language of symbolic dynamics, it is similar to the follower set edge-shift graph of the corresponding shift space, see for example Chapter 3.2 of \cite{LM21}.) We summarize the connection between Penney's ante and this graph in the following lemma.

\begin{lemma} \label{eatpath} Let $G_{v,w}$ be the graph appearing in Figure \ref{1100_1010_graph}. The set of paths in $G_{v,w}$ starting at vertex $(\emptyset, \emptyset)$ and ending at $(v, 100)$ is in bijection with $\Omega_{v,w}(v)$, and the set of paths in $G_{v,w}$ starting at state $(\emptyset, \emptyset)$ and ending at $(10, w)$ is in bijection with $\Omega_{v,w}(w)$. Additionally, the hitting time of the set $\{(v, 100), (10, w)\}$ by the random walk on $G_{v,w}$ which starts at $(\emptyset, \emptyset)$ and assigns probability $p$ to directed edges labeled $1$ and probability $1-p$ to directed edges labeled $0$ has the law of $\min \{\tau_v , \tau_w\}$ under $\p_p$. \end{lemma}

\begin{proof}
Let $\mathcal{Q}(v)$ and $\mathcal{Q}(w)$ denote the set of edge-paths in $G_{v,w}$ starting at $(\emptyset, \emptyset)$ and ending at $(v, 100)$ or $(10, w)$, respectively, and set $\mathcal{Q} = \mathcal{Q}(v) \cup \mathcal{Q}(w)$. Write $\gamma = (e_1, e_2, \ldots, e_t)$ for such a path in $\mathcal{Q}$, and let $\pi$ denote the function from the edge set of $G_{v,w}$ to $\{0,1\}$ which returns the edge label. We claim that the map $\Pi: \mathcal{Q} \to \Omega_{v,w}$ given by 

\begin{equation} \Pi(\gamma) = \pi(e_1) \cdots \pi(e_t) \end{equation}

\begin{figure} 
\includegraphics[width=\textwidth]{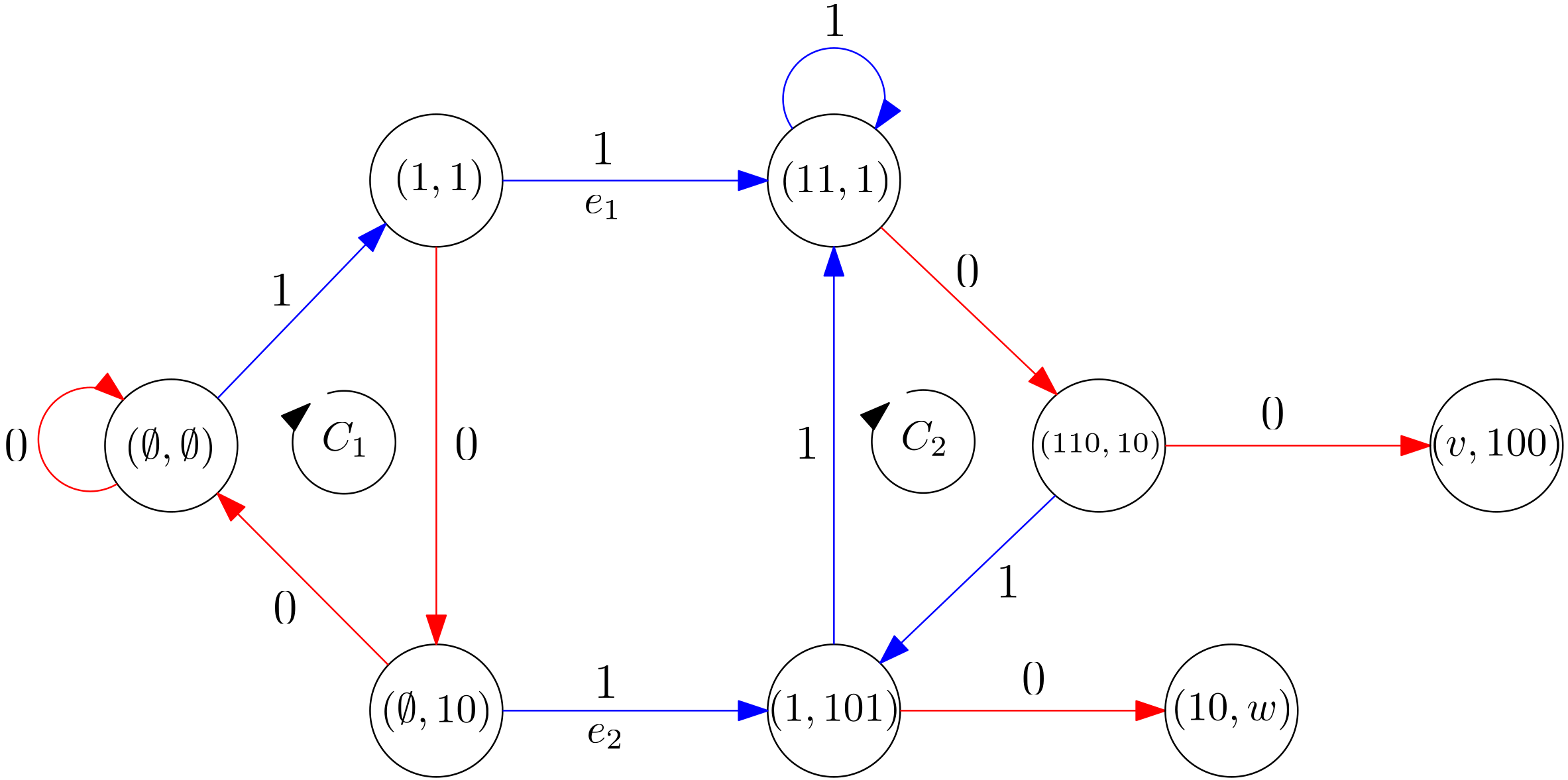} 
\centering
\captionsetup{margin=1cm}
\caption{{\small{The labeled graph $G_{v,w}$ encoding the space $\Omega_{v,w}$, for $v = 1100, w = 1010$. Elements of $\Omega_{v,w}$ are in bijection with labeled paths in $G_{v,w}$ starting at $(\emptyset, \emptyset)$ and ending at either $(v, 100)$ or $(10,w)$, by reading the labels along the edges of the path. (Blue edges are labeled 0, and red edges are labeled 1.)}}}\label{1100_1010_graph}
\end{figure}
\begin{algorithm}
\caption{Construction of $G_{v,w}$}\label{sly_alg}
\begin{algorithmic}[5]
\State {\em{queue}} Stack $\gets [(\emptyset, \emptyset)]$ 
\State {\em{set}} Vertices $ \gets \{\}$ \Comment{Vertices are pairs $(v', w')$ of prefixes of $v, w$, resp.}
\State {\em{set}} Edges $ \gets \{\}$
\While{len(Stack) $> 0$} \Comment{Search the space of possible prefix pairs}
	\State $(v', w') \gets$ pop Stack 
	\For{$a \in \{0,1\}$} \Comment{Append a letter on the right to obtain new prefixes}
		\State $v'' \gets \max \mathcal{O}(v'a, v)$ \Comment{$\max \mathcal{O}$ denotes the maximum length word in $\mathcal{O}$}
		\State $w'' \gets \max \mathcal{O}(w'a, w)$  
		\State add $(v', w') \to (v'', w'')$ to Edges \Comment{This edge is labelled $a$}
		\If{$(v'', w'') \notin$ Vertices}
			\State add $(v'', w'')$ to Vertices
			\If{$v'' \neq v$ {\bf{and}} $w'' \neq w$}: 
				\State add $(v'', w'')$ to Stack
			\EndIf			
		\EndIf
	\EndFor
\EndWhile
\end{algorithmic}
\end{algorithm}

is a bijection with $\Pi(\mathcal{Q}(v)) = \Omega_{v,w}(v)$ and $\Pi(\mathcal{Q}(w)) = \Omega_{v,w}(w)$. To see why, we explain the construction of $G_{v,w}$ via the algorithmic exploration process in Algorithm \ref{sly_alg}. In a nutshell, this algorithm generates all words in $\Omega_{v,w}$ by starting from the empty word, appending one letter at a time, and keeping track of the maximal prefixes of $v, w$ that occur as suffixes of the generated word. 

We prove the result for $\mathcal{Q}(v)$; the result for $w$ is similar. Let $z \in \Omega_{v,w}(v)$. By the algorithmic construction, each vertex in $G_{v,w}$ (except the terminal vertices $(v, 100)$ and $(10, w)$) has exactly two outgoing edges, one for each label in $\{0, 1\}$. So there exists a unique edge-path in $G_{v,w}$ obtained by choosing the edges labeled $z_1, z_2, \ldots$ in succession. Since $z$ has $v$ as a suffix, the resulting path ends at vertex $(v,100)$. Thus $\Pi^{-1}$ is well-defined on $\Omega_{v,w}(v)$, i.e. $\Pi$ is surjective as a map $\mathcal{Q}(v) \to \Omega_{v,w}(v)$. By a direct check of all paths of length 4 in the graph $G_{v,w}$, (or alternatively, because of the condition in Line 12 of Algorithm \ref{sly_alg}), for any $\gamma \in \mathcal{Q}(v)$, $\Pi(\gamma) \in \Omega_{v,w}(v)$. Thus $\Pi$ is well-defined, and clearly injective, so $\Pi$ is the desired bijection. Clearly the pullback of $\p_p$ under $\Pi$ is the (law of the) desired random walk on (edges of) $G_{v,w}$, so the fact that $\Pi^{-1}(\Omega_{v,w}(v)) = \mathcal{Q}(v)$ gives the last conclusion. 

\end{proof}

We are now ready to prove Proposition \ref{prop:quad}. 

\begin{proof} Let $\mathcal{P}(v)$ denote the set of vertex-paths $(\emptyset, \emptyset) \to (v, 100)$ in $G_{v,w}$, and for integers $j, k \geq 0$, let $\mathcal{P}_{j,k}(v)$ denote the set of such paths with $j$ edges labeled 0 (red edges) and $k$ edges labeled 1 (blue edges). We will define a bijection $\psi: \mathcal{P}(v) \to \mathcal{P}(v)$, with the property that for all $j, k \geq 0$, $\psi$ restricted to $\mathcal{P}_{j,k}$ is a bijection $\psi_{j,k}: \mathcal{P}_{j,k}(v) \to \mathcal{P}_{k,j}(v)$. Combining this with Lemma \ref{eatpath} gives the result.

For a path $x \in \mathcal{P}_{j,k}(v)$, let $N_1(x)$ (resp. $N_2(x)$) denote the number of times $x$ visits the vertex $(\emptyset, \emptyset)$ (resp. the vertex $(11, 1)$). $N_1-1$ and $N_2-1$ count the number of times $x$ traverses the cycles $C_1$ and $C_2$, up to a small error. Also, for $1 \leq i \leq N_1(x)$, let $s^1_i(x) \in \{0, 1, 2, \ldots\}$ denote the number of times the self loop at $(\emptyset, \emptyset)$ was consecutively traversed by $x$ immediately after the $i$th visit of $x$ to $(\emptyset, \emptyset)$. Similarly, for $1 \leq i \leq N_2(x)$, let $s^2_i(x)$ denote the number of times the self loop at $(11, 1)$ was consecutively traversed by $x$ immediately after the $i$th visit of $x$ to $(11, 1)$. Note there is no path from any vertex in the 3-cycle $C_2$ to any vertex in the 3-cycle $C_1$, so $x$ traverses exactly one of the two edges $e_1, e_2$ exactly once. Let $e(x) = 1$ or $2$ depending on which edge was traversed. Altogether, any path $x \in \mathcal{P}$ is uniquely described by the tuple

\begin{equation} (N_1(x), N_2(x), (s^1_i(x))_{1 \leq i \leq N_1(x))}, (s^2_i(x))_{1 \leq i \leq N_2(x))}, e(x)), \end{equation}

and each such tuple with $N_1, N_2 \geq 1, s^1_i, s^2_i \geq 0, e \in \{1, 2\}$ determines a unique $x$. Now we define $\psi(x)$ as the path $y \in \mathcal{P}$ given by the tuple

\begin{equation} \label{swapmap} (N_1(y), N_2(y), (s^1_i(y)), (s^2_i(y)), e(y)) = (N_2(x), N_1(x), (s^2_i(x)), (s^1_i(x)), e(x)). \end{equation}

In words, whatever $x$ does in the cycle $C_1$, $y = \psi(x)$ does the same in the cycle $C_2$ and vice versa, and $y$ crosses between the cycles at the same edge that $x$ does. It remains to check that for $x \in \mathcal{P}_{j,k}, y \in \mathcal{P}_{k,j}$. By counting edge labels in $G_{v,w}$ along any path $z \in \mathcal{P}_{r, b}$, one finds

\begin{equation} \label{eq:reds} r = 2 (N_1(z) - 1) + \sum_{i=1}^{N_1(z)} s^1_i(z) + 1\{e(z) = 2\} + N_2(z) + 1, \end{equation}

and 

\begin{equation} \label{eq:blues} b = N_1(z) + 1 + 1\{e(z) = 2\} + \sum_{i=1}^{N_2(z)} s^2_i(z) + 2(N_2(z) - 1). \end{equation}

The terms in Equation \ref{eq:reds} come from: the cycle $C_1$, which is fully traversed $2(N_1(z) - 1)$ times, and contains two edges labeled 0; the self-loop at $(\emptyset, \emptyset)$, which is traversed $\sum s^1_i(z)$ many times; the edge $(1,1) \to (\emptyset, 10)$, which is traversed an additional time (not counted by $N_1(z) - 1$) if $e(z) = 2$; the cycle $C_2$, which is fully traversed $N_2(z)$ times and contains one edge labeled $0$; and the final edge $(110, 10) \to v$, which is traversed exactly once. The counting is similar for Equation \ref{eq:blues}. 

Observe that Equations \ref{eq:reds} and \ref{eq:blues} are exchanged when $N_1(z)$ and $N_2(z)$ are swapped and the sequences $(s^1_i(z))$ and $(s^2_i(z))$ are swapped. Thus, using the definition \ref{swapmap} of $\psi$ and equations \ref{eq:reds} and \ref{eq:blues} for both $(r,b) = (j,k)$ with path $x$ and $(r, b) = (k,j)$ with path $y = \psi(x)$ completes the proof. 

\end{proof} 


\vspace{6pt}
\section{Discussion} \label{open} 

Our work suggests many questions and possible further directions. While our findings describe properties of the win probability for words of different lengths, and for arbitrary $p \in (0,1)$, they say little about the optimal strategy in Penney's ante: that is, given $p$ and a word $w$, which word $v$ optimizes the win probability? Can the existing results on the optimal strategy for the $p = 1/2$ case be extended to arbitrary $p$, or to other underlying Markov chains? Proving Conjecture \ref{baymax} would be a good starting point. Alternatively, it may be more natural to work with different conditions on pairs: for example, rather than requiring $||v||_1 = ||w||_1$, one could restrict to pairs with $|\p_p(v) - \p_p(w)| < \eps$ or $\frac{\p_p(v)}{\p_p(w)} \in (1-\eps, 1+\eps)$ for some $\eps$ depending on $p$ and the lengths of $v$ and $w$. What are the optimal pairings under this type of restriction?

Along the same line, Figure \ref{fig:monoconj}, which plots the volume of pairs where the longer word is favorable in terms of $p$, suggests that in the limit $n \to \infty$ there is a fractal set of $p$ values where the proportion jumps. Specifically, for positive integer $n$, $p \in (0,\frac{1}{2})$ and $\alpha \in (0,1)$, let $a(n, \alpha, p)$ denote the number of pairs $(v,w)$ of length at most $n$ where $v$ is longer than $w$ and Win$(v,w;p) \geq \alpha$, and note that the total number of pairs $(v,w)$ of length at most $n$ where $v$ is longer than $w$ is asymptotically $\frac{1}{3} 4^{n+1}$. 

\begin{question} Does the limit 

\begin{equation} \psi(\alpha, p) = \lim_{n \to \infty} 4^{-n} a(n, \alpha, p) \end{equation}

exist? If so, for what values of $\alpha$ is it a piecewise-continuous, monotone function of $p$? What is the structure of the set of discontinuities? \end{question}

Computer simulations suggest that there is a constant $\alpha_c \approx \frac{1}{3}$ such that $\psi(p)$ is non-increasing for $\alpha < \alpha_c$ and not monotone for $\alpha > \alpha_c$, though we were not able to simulate efficiently enough to properly investigate these questions. 

%

As noted in the introduction, Win$(v,w;p)$ need not be a monotone function of $p$ for fixed $v,w$. In general, what kinds of functions $f: [0,1] \to [0,1]$ can arise as win probability functions, and which cannot? By Theorem \ref{calc1} $f$ is rational, and by Theorem \ref{bread} $f(0), f(1) \in \{0, \frac{1}{2}, 1\}$. What can be said about the derivative of the win probability as a function of $p$? Based on heuristic arguments and computer checks, we believe the following conjecture holds: 

\begin{conjecture} There exists $C > 0$ such that for any positive integer $k$, there exists a pair $(v,w)$ with lengths at most $exp(Ck)$ such that the derivative $\frac{d}{dp}${\normalfont{\text{Win}}}$(v,w;p)$ changes sign at least $k$ times for $p \in (0,1)$. \end{conjecture} 

We also wonder if there is some combinatorial significance to the values of $p$ where this derivative vanishes. Do values of $p$ where $\frac{\partial}{\partial p} \text{Win}(v,w;p) = 0$ correspond to probability spaces with some special or extremal properties, or measures $\p_p$ that solve some variational problem over the space $\Omega_{v,w}$? 

Our results for odd symmetry deal with a majority of cases, but there are still some outliers. The pair $v = 1000, w = 0110$ is an example with odd symmetry, but where property $\mathcal{R}$ does not hold. It is easy to see that there is no length-preserving bijection $f: \Omega_{v, w}(v) \to \Omega_{v,w}(w)$ that transforms $\p_p$ to $\p_{1-p}$: indeed, any length-preserving $f$ must have $f(v) = w$ since these are the only allowable words of length $4$, but $\p_p(v) \neq \p_{1-p}(w)$. The next best thing would be to find a pair of partitions $\mathcal{A}$ of $\Omega_{v,w}(v)$ and $\mathcal{B}$ of $\Omega_{v,w}(w)$ into finite parts, and a bijection $\varphi: \mathcal{A} \to \mathcal{B}$ with $\p_p(A) = \p_{1-p}(\varphi(A))$ for $A \in \mathcal{A}$. We do not know if such a partition exists. 

\begin{problem} Characterize the pairs $(v,w) \notin \mathcal{R}$ that exhibit odd symmetry, and give explicit bijections that explain the symmetry. \end{problem}

The explicit bijections for even and constant symmetry appearing in Section \ref{lotto} suggest many further avenues. In computer-assisted searches, we discovered many pairs that exhibit even and constant symmetry, and which resemble the pairs in $\mathcal{E}$. Specifically, pairs of the form $(ab, \overline{b}a)$ and $(a\overline{a}, b\overline{b})$ (and similar constructions) often have these symmetries, but constructions in the same vein as the maps $\varphi_{a,b}$ do not yield bijections in these cases. Overall, it seems there is more combinatorial structure lurking here.

\begin{problem} Describe a large class of pairs $(v,w) \notin \mathcal{R} \cup \mathcal{E}$ that exhibit even or constant symmetry, and construct explicit bijections to explain the symmetry. \end{problem} 

Here `large' should mean that the number of pairs of length $n$ is $\sim n^{-\frac{1}{2}} 4^{-n}$. (The $\sqrt{n}$ factor arises because of Theorem \ref{bread}, which forces $||v||_1 = ||w||_1$ in the case of even or constant symmetry, and the density of such pairs is roughly $n^{-\frac{1}{2}}$.)

\section*{Acknowledgements} We thank the The Research Experience program (REX) at the University of British Columbia, where this project originated in the Fall of 2022. JR was partially supported by NSERC, ERC Synergy under Grant No. 810115 - DYNASNET and ERC Consolidator Grant 772466 “NOISE".

\bibliographystyle{alpha}
\bibliographystyle{abbrv}
\bibliography{bib}

\end{document}